\documentclass{amsart}
\usepackage{amssymb}
\usepackage{amscd}
\usepackage{verbatim}

\newcommand\nin{\notin}

\newcommand\Id{\operatorname{Id}}
\newcommand\Real{\mathbb{R}}

\newcommand\Cx{\mathbb{C}}

\newcommand\im{\operatorname{Im}}
\newcommand\re{\operatorname{Re}}

\newcommand\Span{\operatorname{span}}

\newcommand\pa{\partial}
\newcommand\Rn{\Real^n}

\newcommand\sphere{\mathbb{S}}
\newcommand\Sn{\sphere^{n-1}}

\newcommand\calE{{\mathcal E}}

\newcommand\calR{{\mathcal R}}

\newcommand\calF{{\mathcal F}}
\newcommand\calG{{\mathcal G}}

\newcommand\Cinf{{\mathcal C}^{\infty}}

\newcommand\Sch{{\mathcal S}}

\newcommand\supp{\operatorname{supp}}

\newcommand\Fr{{\mathcal F}}
\newcommand\Frinv{\Fr^{-1}}
\newcommand\bop{{\mathcal B}}

\newcommand\WF{\operatorname{WF}}

\newcommand\Pt{\tilde P}

\newtheorem{lemma}{Lemma}[section]
\newtheorem{prop}[lemma]{Proposition}
\newtheorem{thm}[lemma]{Theorem}
\newtheorem{cor}[lemma]{Corollary}

\newtheorem*{thm*}{Theorem}
\newtheorem*{prop*}{Proposition}
\newtheorem*{cor*}{Corollary}
\newtheorem*{conj*}{Conjecture}
\numberwithin{equation}{section}
\theoremstyle{remark}
\newtheorem{rem}[lemma]{Remark}
\newtheorem*{rem*}{Remark}
\theoremstyle{definition}

\newtheorem*{Def*}{Definition}

\begin{document}
\renewcommand{\theenumi}{\roman{enumi}}
\renewcommand{\labelenumi}{(\theenumi)}

\title[Fixed energy inverse problem]
{Fixed energy inverse problem for exponentially decreasing potentials}
\author[Gunther Uhlmann and Andras Vasy]{Gunther Uhlmann and Andr\'as Vasy}
\date{November 23, 2001. Addendum added on July 17, 2003.}
\address{Department of Mathematics, University of Washington, 
Seattle, WA}
\email{gunther@math.washington.edu}
\address{Department of Mathematics, Massachusetts Institute of Technology,
Cambridge MA 02139}
\email{andras@math.mit.edu}
\thanks{G.\ U.\ is partially supported by NSF grant \#DMS-00-70488 and
a John Simon Guggenheim fellowship.
A.\ V.\ is partially supported by NSF grant \#DMS-99-70607. Both authors
are grateful for the hospitality of
the Mathematical Sciences Research Institute in Berkeley, CA}

\maketitle

\section{Introduction}
In this paper we show that in two-body scattering the
scattering matrix at a fixed energy
determines real-valued
exponentially decreasing potentials. This result has been proved by Novikov
previously \cite{Novikov:Fixed}, see also
\cite{Novikov-Khenkin:D-bar},
using a $\overline{\pa}$-equation. We present a different
method, which combines a density argument
and real analyticity in part of the complex momentum. The latter has been
noted in \cite{Novikov-Khenkin:D-bar}; here we give a short proof using
contour deformations, similarly to \cite[Section~1.5]{RBMGeo}.
We thus prove:

\begin{thm}\label{thm:main}
Suppose that $n\geq 3$,
$V,V'\in e^{-\gamma_0|w|}L^\infty(\Rn_w;\Real)$ for some $\gamma_0>0$,
and $\lambda>0$. If $S_+(\lambda)=S'_+(\lambda)$, then $V=V'$.
Here $S_+(\lambda)$, resp.\ $S'_+(\lambda)$ are the scattering matrices
of $H=\Delta+V$ and $H'=\Delta+V'$ at energy $\lambda$.
\end{thm}

Theorem~\ref{thm:main} for
compactly supported potentials follows from an analogous result in
\cite{Sylvester-Uhlmann:Global} for
the corresponding Dirichlet-to-Neumann map. See
\cite[Section~12]{Uhlmann:Inverse},
and the references given in these papers for a review of
the relation between the Dirichlet-to-Neumann map and the fixed energy
problem.

The general method follows \cite{Sylvester-Uhlmann:Global}, as discussed
in \cite{RBMGeo}. We thus recall the construction of complex exponential
solutions $u_\rho$, $\rho\in \Cx^n$ of
$(H-\lambda)u_\rho=0$, where $u_\rho(w)=e^{i\rho\cdot w}
(1+v_\rho(w))$, $\rho\cdot\rho=\lambda$, and $v_\rho\to 0$ in an appropriate
sense as $\rho\to\infty$. These solutions exist for $\rho$ outside an
`exceptional set' which is discrete in $z$.
We also show that if we write $\rho=z\nu+\rho_\perp$,
$\nu\in\Sn$, $\rho_\perp\in\Rn$ perpendicular to $\nu$, and $z\in\Cx\setminus
\Real$, then for fixed $\nu$,
$u_\rho$ is analytic in $z$ and real analytic in $\rho_\perp$, hence
extend to be analytic in a neighborhood of $\Real^{n-1}\setminus\{0\}$
in $\Cx^{n-1}_{\rho_\perp}$. The exceptional set is then given by the
zeros of an anlytic function of $z$ and $\rho_\perp$. We caution the reader
that the extension of $u_\rho$ to complex $\rho_\perp$ does not agree
with $u_{z\nu+\rho_\perp}$ where $\rho_\perp$ is allowed to be complex;
indeed $v_\rho$ will merely lie in $e^{\gamma|w|}L^2(\Rn)$ for some $\gamma>0$.

We use this in the inverse problem as follows. Let $u_\rho$, $u'_{\rho'}$
be exponential eigenfunctions of $H$, resp.\ $H'$, as above. Now consider
the pairing
\begin{equation*}
\int_{\Rn} u_\rho (V-V')u'_{\rho'},
\end{equation*}
where $\rho=z\nu+\rho_\perp$, $\rho'=z'\nu+\rho'_\perp$, and $\nu$ is fixed.
If $u_\rho$, $u'_\rho$ are replaced by tempered distributional eigenfunctions
of $H$ and $H'$, then a standard argument shows that $S_+(\lambda)=
S'_+(\lambda)$ implies that the corresponding pairing vanishes.
We employ a density argument to deduce that the pairing also vanishes
for the complex exponentials provided that $|\im z+\im z'|$ is small
and $\rho\cdot\rho=\lambda=\rho'\cdot\rho'$. We then let $\rho,\rho'\to\infty$.
By analyticity, the pairing still vanishes. On the other hand,
$v_\rho, v'_{\rho'}\to 0$, so for $\zeta=\rho-\rho'\in\Real^n$
we deduce that
$\int_{\Rn} e^{i\zeta\cdot w} (V-V')=0$, i.e.\ the Fourier transform
of $V-V'$, hence $V-V'$, vanish. In fact, this step will be
slightly more complicated, because the density argument imposes
restrictions on $\zeta$, and we first deduce vanishing of the Fourier
transform of $V-V'$ in a spherical shell of finite `thickness',
and then use the exponential
decay of $V-V'$ to conclude that it is in fact identically zero.

The authors are grateful to Maciej Zworski for his generous encouragement.

\section{Exponential eigenfunctions}
In this section we recall the construction of
exponential solutions of $(H-\lambda)u=0$ from \cite{Sylvester-Uhlmann:Global}.
First, for $\rho\in\Cx^n$, let
\begin{equation*}
u^0_\rho(w)=e^{i\rho\cdot w}.
\end{equation*}
Thus, $u^0_\rho$ is an `exponential eigenfunction' of $\Delta$, namely
\begin{equation*}
(\Delta-\lambda)u^0_\rho=0,\qquad \rho\cdot\rho=\lambda.
\end{equation*}
We assume everywhere that $n\geq 3$.

For the Hamiltonian $H$, we then seek exponential solutions $u$ of the form
\begin{equation}
u=u_\rho
=e^{i\rho\cdot w}(1+v_\rho),\ \rho\cdot\rho=\lambda,
\ \rho\in \Cx^n,
\end{equation}
where $v_\rho$ is considered a perturbation. In fact, we will have $v_\rho
\in L^2_r(\Rn)$ for all $r<0$ (when $v_\rho$ exists).
Here $L^2_p=L^2_p(\Rn)$ denotes the $L^2(\Rn,\langle w\rangle^{2p}\,dw)$,
$\langle w\rangle^s=(1+|w|^2)^{s/2}$.
Substituting $u$ into $(H-\lambda)u=0$, we obtain
\begin{equation}
(\Delta+2\rho\cdot D_w+V)v_\rho=-V.
\end{equation}
The construction given below works under power-decay assumptions on $V$,
but we state it for exponentially decaying $V$, since $v_\rho$ is real
analytic in the appropriate components of $\rho$ only in that case.
So we assume that
\begin{equation}
V\in e^{-\gamma_0|w|}L^\infty(\Rn),\ \gamma_0>0.
\end{equation}

Thus, we need to construct a right inverse $G(\rho)$ to
\begin{equation*}
P(\rho)=\Delta+2\rho\cdot D_w+V
\end{equation*}
that can be applied to rapidly decreasing functions.
Once this is done,
\begin{equation*}
u_\rho=e^{i\rho\cdot w}(1-G(\rho) V)
\end{equation*}
is the solution to the original problem.
Below we write
\begin{equation}\label{eq:P_0-def}
P_0(\rho)=\Delta+2\rho\cdot D_w.
\end{equation}
Since a right inverse $G_0(\rho)$ of $P_0(\rho)$ can be constructed
explicitly, perturbation theory will give the existence of $G(\rho)$.

Namely, let
\begin{equation*}
G_0(\rho)=\Frinv (|\xi|^2+2\rho\cdot\xi)^{-1} \Fr,
\end{equation*}
so $P_0(\rho)G_0(\rho)=\Id$ e.g.\ on Schwartz functions.
Thus, on the Fourier transform side
$G_0(\rho)$ acts via multiplication by
$(|\xi|^2+2\rho\cdot\xi)^{-1}$ which is in $L^1(\Rn)$.
It is convenient to represent $\rho$ as
\begin{equation*}
\rho=z\nu+\rho_\perp,\ \rho_\perp\in \Rn,\ \nu\in\Sn,
\ z\in\Cx,\ \rho_\perp\cdot\nu=0.
\end{equation*}
We often identify $\Span\{\nu\}^\perp$ with $\Real^{n-1}$.
For $z\nin\Real$, this distribution
is conormal to
\begin{equation}\begin{split}
S(\rho)&=\{\xi\in \Rn:\ |\xi|^2+2\re\rho\cdot\xi=0,
\ \im\rho\cdot\xi=0\}\\
&=\{\xi\in \Rn:\ (\xi+\rho_\perp)^2=\rho_\perp^2,
\ \nu\cdot\xi=0\}.
\end{split}\end{equation}
Note that $S(\rho)$ actually depends only on $\rho_\perp$ and $\nu$,
not on $z$.

Below we write $e^{\gamma\langle w\rangle}L^2(\Rn)=L^2(\Rn;
e^{-2\gamma\langle w\rangle}\,dw)$, and if $m$ is an integer,
\begin{equation*}
e^{\gamma\langle w\rangle}H^m(\Rn)=\{u\in e^{\gamma\langle w\rangle}L^2(\Rn):
\ D^\alpha u\in e^{\gamma\langle w\rangle}L^2(\Rn),\ |\alpha|\leq m\}.
\end{equation*}
The latter is equivalent to $e^{-\gamma\langle w\rangle}u\in H^m(\Rn)$,
hence the notation.

We first recall:

\begin{prop}\cite[Proposition~3.1]{Sylvester-Uhlmann:Global},
\cite[Theorem~1.1]{Weder:Generalized}
$G_0(\rho):L^2_p\to L^2_r$ is bounded for $p>0$, $r<0$, $r\leq p-1$. Moreover,
the norm of $G_0(\rho)$ as a bounded operator between these spaces goes to
$0$ as $|\rho|\to\infty$.
\end{prop}

Our central result is the following proposition.

\begin{prop}\label{prop:analytic-ext}
Suppose that $\gamma>0$ and fix $\nu\in\Sn$.
Then there exists a neighborhood $U$ of $\Real^{n-1}\setminus\{0\}$
in $\Cx^{n-1}$
and an operator
\begin{equation*}
\calG_0(z,\rho_\perp):e^{-\gamma\langle w\rangle}L^2(\Rn)\to
e^{\gamma\langle w\rangle}H^2(\Rn)
\end{equation*}
defined on $(\Cx\setminus\Real)\times U$
such that $\calG_0$ is analytic on $(\Cx\setminus\Real)\times U$
as a bounded operator between these spaces,
and its restriction to $(\Cx\setminus\Real)\times
(\Real^{n-1}\setminus\{0\})$ is $G_0(\rho)$, $\rho=z\nu+\rho_\perp$.
Thus, for $z\in\Cx\setminus\Real$, $\rho_\perp\in\Real^{n-1}\setminus\{0\}$,
the operator $G_0(\rho):e^{-\gamma|w|}L^2\to
e^{\gamma|w|}L^2$ is complex-analytic in $z$, real analytic in $\rho_\perp$.
Moreover, $G_0(\rho)\to 0$ as a bounded operator
on this space as $|\rho|\to\infty$.
\end{prop}

\begin{proof}
We fix some $(z^0,\rho_\perp^0)$, and show that $G_0(\rho)$ extends
to be complex analytic in a neighborhood of this in
$\Cx_z\times\Cx^{n-1}_{\rho_\perp}$. In fact, it is convenient to consider
\begin{equation*}
R_0(\rho)=e^{-i\rho_\perp\cdot w}G_0(\rho)e^{i\rho_\perp\cdot w}.
\end{equation*}
Since the multipliers are holomorphic as maps
\begin{equation*}
e^{\gamma\langle w\rangle}H^m(\Rn)\to e^{\gamma'\langle w\rangle}H^m(\Rn),
\ \gamma<\gamma',
\end{equation*}
for $|\im\rho_\perp|$ sufficiently small, and unitary for $\rho_\perp$ real,
the original statement follows
after we show that $R_0(\rho)$ extends analytically.

We do so by contour deformation on the Fourier transform side.
Let $\xi=(\xi_\parallel,\xi_\perp)$ be the decomposition of $\xi$
according to the decomposition $\Span(\{\nu\})\oplus\Span(\{\nu\})^\perp$
of $\Rn$.
Thus, $\Fr R_0(\rho)\Frinv$ is a multiplication operator by $F^{-1}$
where
\begin{equation}
F(\xi,z,\rho_\perp)=|\xi-\rho_\perp|^2+2(\xi-\rho_\perp)\cdot\rho
=\xi_\parallel^2+2z\xi_\parallel+\xi_\perp^2-\rho_\perp^2.
\end{equation}
Then
\begin{equation*}
\im F=2\im z\,\xi_\parallel,\ \re F=\xi_\parallel^2+2\re z\,\xi_\parallel
+\xi_\perp^2-\rho_\perp^2.
\end{equation*}
Thus the multiplication operator by $F^{-1}$ is singular where $F=0$,
i.e.\ at
\begin{equation*}
\tilde S(\rho)=\{\xi:\ \xi_\parallel=0,\ \xi_\perp^2=\rho_\perp^2\},
\end{equation*}
which is a sphere in the hyperplane $\xi_\parallel=0$.

It is convenient to break up $G_0(\rho)$ into two pieces by introducing
a cutoff $\psi\in\Cinf_c(\Real^n)$ that is identically $1$ near $\tilde
S(\rho^0)$.
For instance, we may take
\begin{equation*}
\psi(\xi)=\phi(\xi_\parallel,\xi_\perp^2)
\end{equation*}
with $\phi\in\Cinf_c(\Real^2)$, identically $1$ near $(0,|\rho_\perp^0|^2)$.
Then
\begin{equation*}
R_0(\rho)=R_0'(\rho)+R_0''(\rho),\ R_0'(\rho)=\Frinv
(|\xi|^2+2\xi\cdot\rho-\rho_\perp^2)^{-1}\psi(\xi)\Fr.
\end{equation*}
Then $R_0''(\rho)$ is a Fourier multiplier by the function $(1-\psi(\xi))
F(\xi,z,\rho_\perp)^{-1}$, which is in fact a symbol of order $-2$,
analytic in $z$ and $\rho_\perp$ for $\im\rho_\perp$, hence $R_0''(\rho)$
is analytic as a map $L^2(\Rn)\to H^2(\Rn)$.

To analyze $R'_0(\rho)$,
it is convenient to introduce polar coordinates in $\xi_\perp$:
$\xi_\perp=r\omega$, $|\omega|=1$, $r\geq 0$. Then
$$
F=\xi_\parallel^2+2z\xi_\parallel+r^2-\rho_\perp^2.
$$
Now, by Fubini's theorem $R_0(\rho)f=\Frinv F^{-1}\Fr f$ can be written as
\begin{equation*}\begin{split}
(R_0'(\rho) f)(w)=(2\pi)^{-n}\int_\Real\int_{\sphere^{n-2}}\int_0^\infty
&e^{i r\omega\cdot w_\perp} e^{i\xi_\parallel w_\parallel}
(\xi_\parallel^2+2z\xi_\parallel+r^2-\rho_\perp^2)^{-1}\\
&\qquad\psi(\xi_\parallel,r\omega)(\Fr f)
(\xi_\parallel,r\omega)\,dr
\,d\omega\,d\xi_\parallel.
\end{split}\end{equation*}
We divide the $\xi_\parallel$ integral into two pieces, corresponding to
$\xi_\parallel\geq 0$ and $\xi_\parallel\leq 0$. In each piece, we then deform
the contour of the $r$ integral in a compact set disjoint from $\supp(1-\psi)$
near $r_0=|\rho_\perp^0|$ in such a way that $\im r^2=2\re r\,\im r$
and $\im z\,\xi_\parallel$ have the same sign on the contour.
Note that the integrand
is analytic in $r$ for $\im r$ small and $\xi_\parallel\neq 0$.

Thus, suppose that $\im z>0$. For $\xi_\parallel> 0$,
we deform the contour $[0,+\infty)_r$
near $r_0$ to a curve $\Gamma_+$
so that $\im r\geq 0$ on $\Gamma_+$ and $r_0$ does {\em not} lie on the
$\Gamma_+$. Now, $F$ never vanishes along $\Gamma_+$,
provided that $\rho_\perp$ is close to $\rho_\perp^0$. Thus, extending $\psi$
to be $1$ on $\Gamma_+\setminus[0,+\infty)$, and using that
$\Fr f$ extends to be analytic in a tube $\{\xi\in\Cx^n:\ |\im\xi|<\gamma\}$,
\begin{equation*}\begin{split}
&(2\pi)^{-n}\int_0^\infty\int_{\sphere^{n-2}}\int_0^\infty
e^{ir\omega\cdot w_\perp} e^{i\xi_\parallel w_\parallel}
(\xi_\parallel^2+2z\xi_\parallel+r^2-\rho_\perp^2)^{-1}\\
&\qquad\qquad\qquad\psi(\xi_\parallel,r\omega)
(\Fr f)(\xi_\parallel,r\omega)\,dr\,d\omega\,d\xi_\parallel\\
&=(2\pi)^{-n}\int_0^\infty\int_{\sphere^{n-2}}\int_{\Gamma_+}
e^{ir\omega\cdot w_\perp} e^{i\xi_\parallel w_\parallel}
(\xi_\parallel^2+2z\xi_\parallel+r^2-\rho_\perp^2)^{-1}\\
&\qquad\qquad\qquad\psi(\xi_\parallel,r\omega)
(\Fr f)(\xi_\parallel,r\omega)\,dr\,d\omega\,d\xi_\parallel,
\end{split}\end{equation*}
and on the right hand side we can allow $\rho_\perp$ to become complex,
proving real analyticity of $R'_0(\rho)$
in $\rho_\perp$, and extending it as an analytic family of
operators $\calR'_0(z,\rho_\perp)$. This argument parallels the
analytic continuation argument of \cite[Chapter~1]{RBMGeo}.
It is now easy to see that $\calR'_0(z,\rho_\perp)$ maps into
$e^{\gamma\langle w\rangle}H^2(\Rn)$; indeed, it maps into
$e^{\gamma\langle w\rangle}C^\infty_{\infty}(\Rn)$, where
$C^\infty_{\infty}(\Rn)$ is the space of smooth functions
which are bounded with all
derivatives.

For $\xi_\parallel<0$ we proceed similarly, deforming the contour
$[0,+\infty)_r$ near $r_0$ to a curve $\Gamma_-$
so that $\im r\geq 0$ on $\Gamma_-$ and $r_0$ does {\em not} lie on the
$\Gamma_-$. Again, we deduce real analyticity in $\rho_\perp$.

The last part follows from the preceeding proposition since $e^{-\gamma|w|}L^2
\subset L^2_p\subset L^2_r\subset e^{\gamma|w|}L^2$.

Instead of the explicit contour deformation, we could have used the partial
Fourier transform in $w_\parallel$, to deduce that
\begin{equation*}
G_0(\rho)f=e^{i\rho_\perp\cdot w}
\Frinv_\parallel (\Delta_\perp+\xi_\parallel^2+2z\xi_\parallel
-\rho_\perp^2)^{-1}\Fr_\parallel e^{-i\rho_\perp\cdot w}f
\end{equation*}
is real analytic in $\rho_\perp$ and analytic in $z$ by inserting step
functions $1=H(\xi_\parallel)+H(-\xi_\parallel)$, and using the analyticity
of
\begin{equation*}
(\Delta-\sigma)^{-1}:e^{-\gamma|w_\perp|}L^2(\Span\{\nu\}^\perp)
\to e^{\gamma|w_\perp|}L^2(\Span\{\nu\}^\perp)
\end{equation*}
in $\sigma$.
\end{proof}

\begin{cor}
Suppose that $\gamma,\gamma_0>0$, $V\in e^{-\gamma_0|w|}L^\infty$.
The operator
\begin{equation*}
V\calG_0(z,\rho_\perp)\in\bop(e^{-\gamma|w|}L^2,e^{-(\gamma_0-\gamma)}L^2)
\end{equation*}
is analytic in $z$ and in $\rho_\perp$ as a bounded operator
between these spaces.
\end{cor}

\begin{cor}
Suppose that $V\in e^{-\gamma_0|w|}L^\infty$ and $\gamma_0>2\gamma$, and
let $U$ be as in Proposition~\ref{prop:analytic-ext}.
Then there exists a set
\begin{equation*}
\calE\subset (\Cx\setminus\Real)_z\times U,
\end{equation*}
which is given by the zeros of an analytic function and whose
intersection with
\begin{equation*}
(\Cx\setminus\Real)_z\times(\Real^{n-1}\setminus\{0\})_{\rho_\perp}
\end{equation*}
is bounded, such that
$(\Id+V\calG_0(z,\rho_\perp))^{-1}$ exists in the complement of $\calE$,
and in a neighborhood of every point where it exists,
$(\Id+V\calG_0(z,\rho_\perp))^{-1}$ is analytic
with values in compact operators on $e^{-\gamma|w|}L^2$.
\end{cor}

\begin{proof}
By the preceeding corollary,
$V\calG_0(\rho):e^{-\gamma|w|}L^2\to e^{-(\gamma_0-\gamma)}L^2$ is analytic
in $z$ and $\rho_\perp$. But the inclusion $e^{-(\gamma_0-\gamma)}L^2
\hookrightarrow e^{-\gamma|w|}L^2$ is compact, so $V\calG_0(z,\rho_\perp)$
is an analytic family of compact operators on $e^{-\gamma|w|}L^2$.
Moreover,
as $|z|\to\infty$ or $|\rho_\perp|\to\infty$, $\rho_\perp$ real,
$V\calG_0(z,\rho_\perp)=VG_0(\rho)\to 0$ in norm.
Thus, the conclusion
follows by analytic Fredholm theory.
\end{proof}

We write
\begin{equation*}
\calG(z,\rho_\perp)=\calG_0(z,\rho_\perp)(\Id+V\calG_0(z,\rho_\perp))^{-1},
\ G(\rho)=G_0(\rho)(\Id+VG_0(\rho))^{-1}.
\end{equation*}

We immediately deduce the following result.

\begin{prop}
Suppose that $V\in e^{-\gamma_0|w|}L^\infty$,
\begin{equation*}
v_{z,\rho_\perp}=
-\calG(z,\rho_\perp)V.
\end{equation*}
Then
\begin{equation*}
((\Cx\setminus\Real)\times U)\setminus\calE
\ni(z,\rho_\perp)\mapsto v_\rho
\end{equation*}
is an analytic function, with values in
$e^{\gamma|w|}L^2$, for any $\gamma>0$.
\end{prop}

\begin{cor}\label{cor:pair-analytic}
Let $\nu\in\Sn$.
Suppose that $V,V'\in e^{-\gamma_0|w|}L^\infty$, and let $\calE$, $\calE'$
be the exceptional sets of these two potentials. Then for
$(z,\rho_\perp)\nin\calE$, $(z',\rho_\perp')\nin\calE'$ the pairing
\begin{equation}\label{eq:pair-88pp}
\int_{\Rn} u_\rho(V-V') u'_{\rho'}
\end{equation}
converges if $|\im z+\im z'|<\gamma_0$,
and is analytic in $z,z',\rho_\perp,
\rho'_\perp$.
\end{cor}

\begin{proof}
We consider the strip $|\im z+\im z'|<\gamma_1<\gamma_0$, $\gamma_1>0$.
Let $\gamma\in(0,(\gamma_0-\gamma_1)/2)$.
Then $1+v_\rho$, $1+v'_{\rho'}$ are analytic
in $(z,z',\rho_\perp,\rho'_\perp)$ with values in $e^{\gamma|w|}L^2$.
Hence,
\begin{equation*}
u_\rho(V-V') u'_{\rho'}=e^{i(\rho+\rho')\cdot w}
(V-V')(1+v_\rho)(1+v'_{\rho'})
\end{equation*}
is analytic in $(z,z',\rho_\perp,\rho'_\perp)$
with values in $L^1(\Rn)$. Integration preserves analyticity and
proves the result.
\end{proof}

\section{Density of generalized eigenfunctions}
In this section we relate tempered distributional eigenfunctions
of $H=\Delta+V$ to its exponential eigenfunctions,
constructed in the previous section.

We first introduce some notation. For $\lambda>0$,
the free incoming Poisson operator is given by
$$
\Pt_{+}(\lambda)g=c
\int_{\Sn} e^{-i\sqrt{\lambda}\, w\cdot\omega}
g\,d\omega_a,\quad g\in\Cinf(\Sn),
\ c=\lambda^{\frac{n-1}4}e^{-\frac{n-1}4\pi i}
(2\pi)^{-\frac{n-1}2}.
$$
The Poisson operator of $H$ is then
$$
P_{+}(\lambda)g=\Pt_{+}(\lambda)g-R(\lambda+i0)((H-\lambda)
\Pt_{+}(\lambda)g)=\Pt_{+}(\lambda)g-R(\lambda+i0)V
\Pt_{+}(\lambda)g.
$$
Note that
for $g\in\Cinf(\Sn)$, $V\Pt_{+}(\lambda)g$ is Schwartz, in fact decays
exponentially, hence $R(\lambda+i0)$ can be applied to it.
For $g\in\Cinf(\Sn)$,
\begin{equation}\label{eq:P_+-form}
P_{+}(\lambda)g=e^{-i\sqrt{\lambda}|w|}g_-+e^{i\sqrt{\lambda}|w|}g_++L^2(\Rn),
\ g_+,g_-\in\Cinf(\Sn),\ g_-=g.
\end{equation}
For such $g$, $P_+(\lambda)g$ is characterized by the property that it is
the unique solution $u$ of $(H-\lambda)u=0$ which is of the form
\eqref{eq:P_+-form}.
The scattering matrix is then the operator
\begin{equation*}
S_+(\lambda):\Cinf(\Sn)\to\Cinf(\Sn),\ S_+(\lambda)g_-=g_+.
\end{equation*}

There is also an incoming Poisson operator $P_-(\lambda)$ which is
characterized by the fact that for $g\in\Cinf(\Sn)$,
$P_-(\lambda)g$ is the unique solution $u$ of
$(H-\lambda)u=0$ of the form
\begin{equation}\label{eq:P_--form}
P_{-}(\lambda)g=e^{-i\sqrt{\lambda}|w|}g_-+e^{i\sqrt{\lambda}|w|}g_++L^2(\Rn),
\ g_+,g_-\in\Cinf(\Sn),\ g_+=g.
\end{equation}
In particular, for $g\in\Cinf(\Sn)$,
\begin{equation}\label{eq:cx-conjugate}
\overline{P_-(\lambda)g}=P_+(\lambda)\overline{g}.
\end{equation}

The S-matrix is related to the Poisson operator via the following boundary
pairing.

\begin{prop}\label{prop:boundary-pair}\cite[Lemma~2.2]{RBMGeo}
Suppose that $\lambda>0$, and $V\in e^{-\gamma_0|w|}L^\infty$, $\gamma_0>0$.
Suppose that $(H-\lambda)u_+\in L^2_s$, $(H-\lambda)u_-\in L^2_s$, $s>1/2$,
and
\begin{equation*}\begin{split}
&u_+=e^{-i\sqrt{\lambda}|w|}g_{+-}+e^{i\sqrt{\lambda}|w|}g_{++}+L^2,\\
&u_-=e^{-i\sqrt{\lambda}|w|}g_{--}+e^{i\sqrt{\lambda}|w|}g_{-+}+L^2,
\end{split}\end{equation*}
$g_{\pm\pm}\in\Cinf(\Sn)$.
Then
\begin{equation}
\langle u_+,(H-\lambda)u_-\rangle
-\langle(H-\lambda) u_+,u_-\rangle=
2i\sqrt{\lambda}(\langle g_{++},g_{-+}\rangle-
\langle g_{+-},g_{--}\rangle).
\end{equation}
\end{prop}

\begin{rem}
This is stated for $V\in\Cinf_c(\Rn)$ in \cite{RBMGeo}. However, if $u_+, u_-$
are
as above, then $Vu_+ u_-\in e^{-\gamma|w|}L^1$ for $\gamma<\gamma_0$, hence
the conclusion is equivalent to the corresponding statement with
$H-\lambda$ replaced by $\Delta-\lambda$.
\end{rem}

Let $R(\lambda')=(H-\lambda')^{-1}$ for $\lambda'\in\Cx\setminus\Real$.
Let $f\in\Sch(\Rn)$, $g\in\Cinf(\Sn)$, and apply this proposition with
\begin{equation*}
u_-=R(\lambda-i0)f=e^{-i\sqrt{\lambda}|w|}g_{--}+L^2,\ u_+=P_+(\lambda)g.
\end{equation*}
We deduce that
\begin{equation}\label{eq:dense-pairing}
\langle u_+,f\rangle=-2i\sqrt{\lambda}\langle g,g_{--}\rangle.
\end{equation}

Our density result is the following.

\begin{prop}\label{prop:density}
Suppose that $V\in e^{-\gamma_0|w|} L^\infty$, and let
$0<\gamma<\gamma'<\gamma_0$. Then the set
$$
\calF=\{P_{+}(\lambda)g_+:\ g_+\in\Cinf(\Sn)\}
$$
is dense in the nullspace of $H-\lambda$ on $e^{\gamma|w|}L^2$
in the topology of $e^{\gamma'|w|}L^2$.
\end{prop}

\begin{proof}
Suppose that $f\in e^{-\gamma'|w|}L^2$ is orthogonal to $\calF$.
Let $u_-=R(\lambda-i0)f$. By \eqref{eq:dense-pairing}, for all $g\in
\Cinf(\Sn)$, $\langle g,g_{--}\rangle=0$ since
$\langle f,P_{+}(\lambda)g\rangle$ vanishes by assumption.
But $u_-=R_0(\lambda-i0)f'$, $f'=f-VR_0(\lambda-i0)f
\in e^{-\gamma|w|}L^2$. Thus, $\Fr u_-$ is the product of an analytic function,
namely $\Fr f$,
and $(|\xi|^2-(\lambda-i0))^{-1}$. Thus, $u_-\in L^2$ implies
that $\Fr f$ vanishes on the sphere $|\xi|=\sqrt{\lambda}$. Hence
$\Fr f=(\xi^2-\lambda)\phi$, with $\phi$ analytic in the strip
$|\im\xi|<\gamma'$. Thus, $u_-\in e^{-\gamma|w|}L^2$ for $\gamma<\gamma'$.
Thus for $v\in e^{\gamma|w|}L^2$ with $(H-\lambda)v=0$,
$$
\langle f, v\rangle=\langle (H-\lambda)u_-,v\rangle
=\langle u_-,(H-\lambda)v\rangle=0,
$$
i.e.\ $f$ ia orthogonal to the nullspace of $H-\lambda$ on
$e^{\gamma|w|}L^2$. Thus, $\calF$ is dense in this nullspace.
\end{proof}

Our approach to the inverse problem relies on relating the S-matrices
to the pairing \eqref{eq:pair-88pp}. Thus,
we consider two many-body operators $H$ and $H'$ induced
by potentials $V$ and $V'$ respectively, and show that the equality
of the S-matrices at a fixed energy $\lambda$ implies the vanishing of
an analogous pairing. For this we use the following consequence
of Proposition~\ref{prop:boundary-pair} applied with $\Delta$ in place
of $H$.

\begin{prop}
Suppose that $\lambda>0$.
Let $u_+=P_{+}(\lambda)g_+$, $u_-=P'_{-}(\lambda)g_-$. Then
\begin{equation}
\langle u_+,(\Delta-\lambda)u_-\rangle
-\langle(\Delta-\lambda) u_+,u_-\rangle=
2i\lambda(\langle S_{+}(\lambda)g_+,g_-\rangle-
\langle g_+,S'_{-}(\lambda)g_-\rangle).
\end{equation}
\end{prop}

\begin{cor}\label{cor:pair-8}
Suppose $\lambda>0$, $S_{+}(\lambda)=S'_{+}(\lambda)$.
Let $u_+=P_{+}(\lambda)g_+$,
$u_-=P_{-}(\lambda)g_-$. Then
\begin{equation}\label{eq:pair-7}
\int_{\Rn} (V-V')u_+ \overline{u_-}=0.
\end{equation}
Similarly, if $u_+=P_{+}(\lambda)g_+$,
$u_-=P_{+}(\lambda)g_-$, then.
\begin{equation}\label{eq:pair-8}
\int_{\Rn} (V-V')u_+ u_-=0.
\end{equation}
\end{cor}

\begin{proof}
\eqref{eq:pair-7} follows from the preceeding proposition since
$S'_-(\lambda)^*=S'_+(\lambda)$. Then \eqref{eq:pair-8} follows from
\eqref{eq:pair-7} by applying the latter with $g_-$ replaced by
$\overline{g_-}$ and using \eqref{eq:cx-conjugate}.
\end{proof}

\section{Inverse results: Proof of Theorem~\ref{thm:main}}
Let $\lambda>0$, and suppose that
\begin{equation*}
V,V'\in e^{-\gamma_0|w|}L^\infty,\ \gamma_0>0.
\end{equation*}
Fix $\zeta\in\Rn$ such that $|\zeta|>2\sqrt{\lambda}$,
and let $\nu\in\Sn$ be orthogonal to $\zeta$, and let
$\mu\in\Sn$ orthogonal to both $\zeta$ and $\nu$. For $t$ real,
$t>\sqrt{\frac{1}{4}|\zeta|^2-\lambda}$, let
\begin{equation}\begin{split}\label{eq:rho-param}
&\rho=\rho(t)=\frac{\zeta}{2}+(t^2-\frac{1}{4}|\zeta|^2+\lambda)^{1/2}\mu+it\nu,\\
&\rho'=\rho'(t)=\frac{\zeta}{2}-(t^2-\frac{1}{4}|\zeta|^2+\lambda)^{1/2}\mu-it\nu,
\end{split}\end{equation}
so $\rho\cdot\rho=\lambda=\rho'\cdot\rho'$.
By Corollary~\ref{cor:pair-analytic}, the integral
\begin{equation}\label{eq:pair-88}
\int_{\Rn} u_\rho(V-V') u'_{\rho'}
\end{equation}
converges for all $t$, and is meromorphic in $t$ in a neighborhood of
$(\sqrt{\frac{1}{4}|\zeta|^2-\lambda},+\infty)$.
We use a density argument, Proposition~\ref{prop:density},
and Corollary~\ref{cor:pair-8}
to show that this integral actually vanishes if $S_{+}(\lambda)
=S'_{+}(\lambda)$ and
\begin{equation}\label{eq:zeta-bds}
2\sqrt{\lambda}<|\zeta|<\sqrt{4\lambda+\gamma_0^2}.
\end{equation}

Indeed,
for $\sqrt{\frac{1}{4}|\zeta|^2-\lambda}<t<\gamma<\gamma_0/2$, $u'_{\rho'}$
can be approximated by $P_{-}(\lambda) g_-$
in $e^{\gamma|w|}L^2$ due to Proposition~\ref{prop:density}. Similarly,
$u_\rho$
can be approximated by $P_{+}(\lambda) g_+$
in $e^{\gamma|w|}L^2$.
On the other
hand, $V-V'$ lies in $e^{-\gamma_0|w|}L^2$.
Hence the product can
be approximated in $L^1$ by a product which takes the form
of the integrand of \eqref{eq:pair-8}. The equality of the S-matrices
implies that \eqref{eq:pair-8} vanishes, hence so does \eqref{eq:pair-88},
i.e.\ we deduce the following result.

\begin{prop}
Suppose that $\lambda>0$, $V,V'\in e^{-\gamma_0|w|}L^\infty$,
$S_{+}(\lambda)=S'_{+}(\lambda)$. Then for $\zeta$ satisfying
\eqref{eq:zeta-bds}, $\rho$, $\rho'$ given by \eqref{eq:rho-param}
with $(z,\rho_\perp)\nin\calE$, $(z,\rho_\perp)\nin\calE'$,
\begin{equation}\label{eq:pair-88p}
\int_{\Rn} u_\rho (V-V')u'_{\rho'}=0
\end{equation}
for $\sqrt{\frac{1}{4}|\zeta|^2-\lambda}<t<\gamma_0/2$.
\end{prop}

The pairing in \eqref{eq:pair-88p} is meromorphic in $t$ with $\re t>
\sqrt{\frac{1}{4}|\zeta|^2-\lambda}$ and $|\im t|$ sufficiently small.
It vanishes on an interval inside this domain by the proposition. Thus,
\eqref{eq:pair-88p} holds for all $t>
\sqrt{\frac{1}{4}|\zeta|^2-\lambda}$.
Then as $t\to\infty$, the integral on the left hand side of
\eqref{eq:pair-88p} converges to
\begin{equation}\label{eq:pair-99}
\int_{\Rn} (V-V') u^0_\rho u^0_{\rho'}=\int_{\Rn}
(V-V')\,e^{i\zeta\cdot w}\,dw
\end{equation}
since $v_\rho\to 0$, $v'_{\rho}\to 0$ in $e^{\gamma|w|}L^2(\Rn)$ for
any $\gamma>0$ and $V-V'\in e^{-\gamma_0|w|}L^\infty(\Rn)$ with
$\gamma_0>0$.
But this is the Fourier transform of $V-V'$,
evaluated at $\zeta$.
Hence
the vanishing of \eqref{eq:pair-99} shows that the Fourier transform of
$V-V'$ vanishes on the shell \eqref{eq:zeta-bds}. Since this Fourier
transform is real analytic, as $V-V'\in e^{-\gamma_0|w|}L^\infty$, we deduce
that it vanishes everywhere, hence $V=V'$. This completes the proof
of Theorem~\ref{thm:main}.

\appendix

\section{Addendum, July 17, 2003}

The main result of this paper was a new proof of the 
fixed energy inverse result for exponentially decaying potentials in
potential scattering, which was first obtained by
Novikov \cite{Novikov:Fixed}, see also
the work of Novikov and Khenkin \cite{Novikov-Khenkin:D-bar}.
We used two main techniques:
real analyticity in part of the complex momentum (Proposition~2.2)
and a density argument.
The former had been noted in \cite{Novikov-Khenkin:D-bar}, as we pointed
out in the paper. We provided a simple proof by showing that the
appropriate quantities have an analytic continuation (over which we
needed no control) --
this is similar to the analytic
continuation arguments for the resolvent of the Laplacian as can be found,
for instance, in
Melrose's book \cite[pp.\ 7-9]{RBMGeo}.

What we did not realize that after Novikov, but several years before
we started working on this project, Eskin and Ralston wrote a paper
\cite{Eskin-Ralston:Magnetic} in
which they not only extended Novikov's result to the magnetic setting,
but also simplified his proof in the potential setting.
Unfortunately we were unaware of this aspect of their paper, as well
as of their proofs. It turns out
that the contour deformation proof we gave for real analyticity is
substantially identical
to that of Eskin and Ralston \cite[pp.\ 204-205]{Eskin-Ralston:Magnetic},
although they also provide more control over the analytic continuation,
and use it elsewhere in their paper. We very much regret to have
omitted the reference to \cite{Eskin-Ralston:Magnetic}.

The rest of the proof of the inverse result is different; ours
relies on a density argument, while that of Eskin and Ralston
on taking a limit to real
frequencies.

\def\cprime{$'$}


\begin{thebibliography}{1}

\bibitem{Eskin-Ralston:Magnetic}
G.~Eskin and J.~Ralston.
\newblock Inverse scattering problem for the {S}chr{\"o}dinger equation with
  magnetic potential at a fixed energy.
\newblock {\em Commun. Math. Phys.}, 173:199--224, 1995.

\bibitem{RBMGeo}
R.~B. Melrose.
\newblock {\em Geometric scattering theory}.
\newblock Cambridge University Press, 1995.

\bibitem{Novikov-Khenkin:D-bar}
R.~G. Novikov and G.~M. Khenkin.
\newblock The $\overline\partial$-equation in the multidimensional inverse
  scattering problem.
\newblock {\em Uspekhi Mat. Nauk}, 42(3(255)):93--152, 255, 1987.

\bibitem{Novikov:Fixed}
Roman~G. Novikov.
\newblock The inverse scattering problem at fixed energy for the
  three-dimensional {S}chr\"odinger equation with an exponentially decreasing
  potential.
\newblock {\em Comm. Math. Phys.}, 161(3):569--595, 1994.

\bibitem{Sylvester-Uhlmann:Global}
J.~Sylvester and G.~Uhlmann.
\newblock A global uniqueness theorem for an inverse boundary value problem.
\newblock {\em Ann. of Math.}, 125:153--169, 1987.

\bibitem{Uhlmann:Inverse}
Gunther Uhlmann.
\newblock Inverse boundary value problems and applications.
\newblock {\em Ast\'erisque}, (207):6, 153--211, 1992.
\newblock M\'ethodes semi-classiques, Vol.\ 1 (Nantes, 1991).

\bibitem{Weder:Generalized}
Ricardo Weder.
\newblock Generalized limiting absorption method and multidimensional inverse
  scattering theory.
\newblock {\em Math. Methods Appl. Sci.}, 14(7):509--524, 1991.

\end{thebibliography}
\end{document}